\documentclass[11pt]{amsart}
\usepackage[margin=1in]{geometry}
\usepackage{amssymb,booktabs,mathrsfs}
\usepackage[hidelinks]{hyperref}
\newtheorem{theorem}{Theorem}[section]
\newtheorem{lemma}[theorem]{Lemma}
\newtheorem{proposition}[theorem]{Proposition}
\newtheorem{corollary}[theorem]{Corollary}
\theoremstyle{remark}

\newcommand{\g}{\mathfrak g}
\newcommand{\gr}{\operatorname{gr}}
\newcommand{\Gr}{\operatorname{Gr}}
\newcommand{\cA}{\mathcal A}
\newcommand{\R}{\mathbb R}

\newcommand{\Z}{\mathbb Z}

\title[Quasi-isometric lattices]{Quasi-isometric lattices that are not bi-Lipschitz equivalent}
\author{Dmitri Burago}
\date{November 2026. Revised draft for review.}

\begin{document}
\begin{abstract}
We exhibit a simply connected nilpotent Lie group of dimension eight
in which a cocompact lattice is not bi-Lipschitz equivalent to any of
its proper finite-index subgroups. An exact rational computation
obstructs dilations of the filtered cohomology ring. To connect the
filtration with the geometry, we prove that a top-degree group cocycle
of subcritical growth represents zero. The proof smooths the cocycle
at large scales, and the weighted derivative estimates recover the
homogeneous dimension. We then use the cohomological induction
associated with an extreme first-layer map of a quasi-isometry. This construction was proposed by Sergei Ivanov and myself, and we used two AIs to complete computations. 
\end{abstract}
\maketitle

\section{Introduction}

Cocompact lattices in the same Lie group are quasi-isometric. The same
question appears in \cite[Problem~2]{BK} and \cite[Question~22]{MFO}. Problem
13 of \cite{BCLQ} asks whether they are necessarily bi-Lipschitz
equivalent. We give a negative answer using a lattice and a subgroup of
index two. This construction was proposed by Sergei Ivanov and myself, and we 
used Claude Max and Astra 6 pro to complete computations. The problem had been posed by
Curt McMullen in our conversation with Bruce Kleiner and me. Note that such examples cannot exist in Hyperbolic groups and for lattices with the same co-volume (Hall's 	``Marriage lemma''), so we needed a nilpotent group without volume dilations.

\begin{theorem}\label{thm:main}
Let $G$ be the simply connected nilpotent Lie group whose Lie algebra is
specified in Section~\ref{sec:example}. For every cocompact lattice
$\Gamma<G$ and every proper finite-index subgroup $\Gamma'<\Gamma$,
the lattices $\Gamma$ and $\Gamma'$ are not bi-Lipschitz equivalent.
\end{theorem}

The Lie algebra is rational, so lattices exist by Mal'cev's theorem \cite{Malcev}.
The abelianization of each lattice has free rank two. Composing
a projection onto a free cyclic quotient with reduction modulo two
gives a surjection to $\Z/2$, whose kernel provides the index-two
example. All word metrics on a fixed finitely generated group are
bi-Lipschitz equivalent, so the choice of word metrics is immaterial.

A bijection onto a subgroup of index $k$ fixes the volume scaling of
every blow-down. In the present example the associated Carnot algebra
has a small graded automorphism group, so this fixes the dilation
parameter of every tangent map of a blow-down. We use the induction
construction of Gotfredsen and Kyed \cite{GK} with an ergodic measure
selected by Mukherjee, Samanta and Thandar \cite{MST}. A bound on the
growth of cocycles lets us compute the associated graded of this
induced cohomology map. The resulting dilation violates the rational
cohomology calculation unless $k=1$.

The filtration is important here. In a cohomology ring one can always
multiply the degree-$j$ component by $s^j$, for any $s\ne0$. The
obstruction below prescribes the action induced by the Carnot
dilation on each weight of the associated graded cohomology.

\section{The Lie algebra}\label{sec:example}

Let $\g$ have basis $e_0,\dots,e_7$ with \emph{weights} $(1,1,2,3,3,4,4,5)$, adapted to the lower central series, and the nonzero brackets
\begin{align*}
&[e_0,e_1]=e_2,\quad [e_0,e_2]=e_3,\quad [e_0,e_3]=e_5,\quad [e_0,e_4]=e_6,\quad [e_0,e_5]=e_7,\quad [e_0,e_6]=-\tfrac{29}{33}e_7,\\
&[e_1,e_2]=e_4,\quad [e_1,e_3]=e_6,\quad [e_1,e_4]=-\tfrac19e_5+\tfrac23e_6+\boxed{\tfrac{35}{33}e_7},\quad [e_1,e_5]=-\tfrac{23}{11}e_7,\\
&[e_1,e_6]=-\tfrac{103}{33}e_7,\quad [e_2,e_3]=\tfrac{40}{33}e_7,\quad [e_2,e_4]=\tfrac{80}{33}e_7 .
\end{align*}
Every bracket raises weight by \emph{exactly} the sum of the weights, with one exception: the boxed term (weights $1+3=4$, but $e_7$ has weight $5$). Deleting it gives the associated graded (Carnot) algebra $\gr\g$. The structure constants are rational, so $G$ has lattices by Mal'cev's theorem \cite{Malcev}. The homogeneous dimension is $Q=1+1+2+3+3+4+4+5=23$, and $\dim G=8$.

\begin{proposition}\label{prop:facts}
The displayed brackets satisfy the Jacobi identity. The derivation
algebra of $\g$ has dimension $12$ and every derivation is traceless.
The graded automorphisms of $\gr\g$ are
\[
 \delta_t\quad\text{and}\quad\delta_t\tau \qquad(t\in\R^*),
\]
where $\delta_t e_i=t^{w_i}e_i$ and $\tau$ is an involution of
determinant one. Its first-layer matrix is
$\frac15\begin{pmatrix}-7&-4\\6&7\end{pmatrix}$.
\end{proposition}

\begin{proof}
For Jacobi, substitute the brackets on all triples of basis vectors.
Equivalently, the resulting Chevalley--Eilenberg differential squares
to zero on the generators. For derivations, use the $64$ matrix
entries of $D$ as variables and impose
$D[e_i,e_j]=[De_i,e_j]+[e_i,De_j]$ for all $i<j$.
The rational equation matrix has rank $52$, and the trace row belongs
to its row span. Since the automorphism group is real algebraic, it
has finitely many real connected components. Its absolute determinant
character has zero derivative by the trace computation and therefore
has image $\{1\}$.

For the graded automorphisms, write the first-layer matrix as
$\begin{pmatrix}p&q\\r&s\end{pmatrix}$. The higher images are obtained
successively from the expressions
\[
 e_2=[e_0,e_1],\quad e_3=[e_0,e_2],\quad e_4=[e_1,e_2],\quad
 e_5=[e_0,e_3],\quad e_6=[e_0,e_4],\quad e_7=[e_0,e_5].
\]
Impose bracket preservation on all basis pairs. Divide the resulting
nonzero equations by their common factor $ps-qr$ and append
$z(ps-qr)-1=0$. The branch $s=0$ has Gr\"obner basis $\{1\}$.
Otherwise scale the first-layer matrix so that $s=1$, using a
dilation. In the order $z,p,q,r$, the lexicographic basis is
\[
 z+\tfrac{259}{75}r-1,\qquad
 p+\tfrac73r-1,\qquad q+\tfrac23r,\qquad r^2-\tfrac67r.
\]
The two normalized solutions give the identity and a scalar
multiple of the displayed first-layer matrix. Substitution of the
latter matrix in the bracket expressions gives $\tau^2=1$ and
$\det\tau=1$ on the full algebra. This proves the classification.
\end{proof}

\section{The filtered cohomology obstruction}\label{sec:algebra}

Write $e^i$ for the dual basis and assign $e^i$ weight $w_i$.
Monomial weights add under the wedge product. Let $F_p$ be the
span of the monomials of weight at most $p$. The differential
preserves $F_p$ and has weight-preserving part $d_0$, the
differential for $\gr\g$. Put $W_pH^j(\g)$ equal to the image
of the cocycles in $F_p\cap\Lambda^j\g^*$. The multiplicative
spectral sequence of this filtration starts at
$E_1=H^*(\gr\g)$ and ends at $\Gr^W H^*(\g)$.

\begin{proposition}\label{prop:alg}
Both $\g$ and $\gr\g$ have Betti numbers $(1,2,4,8,10,8,4,2,1)$.
There is a canonical identification of weight-graded rings
$\Gr^W H^*(\g)=H^*(\gr\g)$. The cup product on weights $p$ and
$Q-p$ induces a perfect pairing between
$\Gr^W_p H^j$ and $\Gr^W_{Q-p}H^{8-j}$, where $Q=23$, and
$W_pH^j\smile W_{Q-p-1}H^{8-j}=0$.
\end{proposition}

\begin{proof}
The Betti numbers follow by rational elimination using
$d e^k=-\sum_{i<j}e^k([e_i,e_j])e^i\wedge e^j$ and the Leibniz
rule. The dimensions on $E_1$ and $E_\infty$ agree, so every
spectral-sequence differential vanishes. More explicitly, leading
parts of actual cocycles are contained in $\ker d_0$, and
$\operatorname{im}d_0$ is contained in the leading parts of actual
coboundaries. The equality of the Betti numbers implies equality of
the ranks of $d$ and $d_0$ in each degree. These two inclusions are
therefore equalities. The canonical identification sends a leading
part to its $d_0$-cohomology class. In particular, all graded
automorphisms in Proposition~\ref{prop:facts} act on the associated
graded of the cohomology filtration.

Poincar\'e duality on the associated graded algebra follows from
nilpotent Lie algebra duality or from a compact quotient of the
rational graded group and Nomizu's theorem. The top class has weight
$Q$, so the pairing respects complementary weights and is perfect
there. A product in top degree of total weight below $Q$ vanishes
already on cochains.
\end{proof}

\begin{proposition}\label{prop:rigid}
Suppose a degree-preserving derivation $D$ of $H^*(\g)$ preserves
$W$ and acts on $\Gr^W_pH^j$ as multiplication by $\lambda p$
for every $p,j$. Then $\lambda=0$.
\end{proposition}

\begin{proof}
Choose a basis of cohomology by computing cycles in order of
increasing weight and retaining each new independent class modulo
all coboundaries. Use $h_{j,i}$ for the resulting basis vectors
and $p_{j,i}$ for their weights. Write
\[
 D h_{j,i}=\lambda p_{j,i}h_{j,i}
       +\sum_{a:\ p_{j,a}<p_{j,i}}d^j_{ai}h_{j,a}.
\]
Reduce the wedges of the chosen cocycles modulo coboundaries in
these bases. For every pair of basis classes $u,v$, impose
$D(uv)-(Du)v-u(Dv)=0$ coefficient by coefficient. This is a
homogeneous rational system on $79$ weight-lowering entries and
$\lambda$. Its rank is $51$ and the row $\lambda=0$ reduces to
zero against its row echelon basis. The same check for $\gr\g$ has
rank $50$ and leaves $\lambda$ free, as expected from its dilations.
The companion verifier \texttt{bilip\_lattices\_exact\_check.js}
constructs the differential and products directly from the brackets
and performs these reductions using rational arithmetic.
\end{proof}

\begin{corollary}\label{cor:rigid}
If a filtered graded-ring automorphism $\Phi$ of $H^*(\g)$ has
associated graded equal to $(\delta_t g)^*$, with $t>0$ and
$g\in\{1,\delta_{-1},\tau,\delta_{-1}\tau\}$, then $t=1$.
\end{corollary}

\begin{proof}
The pairs $(\Phi,s)$ with $s\in\R^*$ and
$\Gr\Phi|_{\Gr_p}=s^p\,\mathrm{Id}$ form a real algebraic group.
The derivative of its character $(\Phi,s)\mapsto s$ vanishes
by Proposition~\ref{prop:rigid}. Its image is finite because a
real algebraic group has finitely many connected components,
so that image lies in $\{\pm1\}$. The four listed automorphisms
commute with the dilations and square to the identity. Applying
the character to $(\Phi^2,t^2)$ gives $t=1$.
\end{proof}

\section{Growth of cocycles}\label{sec:growth}

Use the Mal'cev coordinates of the adapted basis and let
$D_R e_i=R^{w_i}e_i$. Write $\sigma_R(x)=D_R^{-1}\log x$,
regarding its values as coordinates in the Carnot group
$G_\infty$ with Lie algebra $\gr\g$. For $R\ge1$, the
transported group law in these coordinates has structure constants
$c_{ij}^{\,l}R^{w_i+w_j-w_l}$. It converges to the Carnot group law
with all derivatives on compact sets, by the finite
Baker--Campbell--Hausdorff formula. The rescaled distances also
converge on compact sets by Pansu's theorem \cite{Pan83,MST}.
In particular, a fixed compact coordinate set corresponds
to points at distance $O(R)$ from the identity, and balls
of radius $O(R)$ rescale into a fixed compact set. We use any
fixed left-invariant Riemannian metric on $G$ for these
comparisons. On $\Gamma$ it is quasi-isometric to the word
metric via the inclusion, so the cocycle growth conditions are
equivalent for the two choices.

Use homogeneous invariant group cochains on $\Gamma$.
We say a cochain has growth $O(R^p)$ or $o(R^p)$ if the
supremum of its absolute value on tuples of diameter at most
$R$ has that bound. A left-invariant form $\omega$ on $G$
gives a cochain
\[
 \xi_\omega(x_0,\ldots,x_j)=
           \int_{\sigma(x_0,\ldots,x_j)}\omega .
\]
Here the simplex $\sigma$ is defined recursively by coning
the opposite face from $x_0$ along exponential segments.
Faces and left translations commute with this construction,
so Stokes' theorem gives a cochain map. On $\Gamma$ this
realizes the Nomizu isomorphism \cite{Nom}. Indeed, the
coned simplices give an equivariant chain map from the
homogeneous bar resolution to the singular chains of the
contractible universal cover, lifting the identity in
degree zero. Integration of forms gives the usual comparison
with de Rham cohomology.

\begin{lemma}\label{lem:integration}
For $\omega\in F_p$, the cochain $\xi_\omega$ has growth
$O(R^p)$. If $\sigma_R(x_i^R)\to z_i$ for each $i$, then
\[
 R^{-p}\xi_\omega(x_0^R,\ldots,x_j^R)
 \longrightarrow
 \int_{\sigma_\infty(z_0,\ldots,z_j)}\omega^{(p)},
\]
where $\omega^{(p)}$ is the weight-$p$ part viewed as a
left-invariant form on $G_\infty$. The convergence of the
rescaled integrals to their Carnot analogues is uniform on
compact sets of rescaled configurations.
\end{lemma}

\begin{proof}
In the transported group law, the pullback of $R^{-p}\omega$
has invariant coefficients multiplied by $R^{q-p}$ in each
weight $q$. These forms converge smoothly on compact sets to
the indicated form. The recursively coned simplices can be
parametrized by a cube, coning with one parameter at a time.
In coordinates, each step consists of multiplication,
inversion and multiplication of the logarithm by a scalar
between zero and one. The parametrizations and their first
derivatives converge uniformly on compact sets of input
tuples. Integrating the convergent pullbacks gives the formula.
For the growth estimate, left translate the first vertex to
the identity and apply boundedness on a compact set of
rescaled tuples.
\end{proof}

\begin{lemma}[Top-degree growth]\label{lem:top}
If a $\Gamma$-invariant top-degree cocycle $\zeta$ has growth
$o(R^Q)$, it represents zero in $H^8(\Gamma;\R)$.
\end{lemma}

\begin{proof}
Replace $\zeta$ by its alternation, which preserves its class
and its growth bound. Normalize Haar measure $dy$ to be
Lebesgue measure in exponential coordinates. Choose a
nonnegative smooth compactly supported function $\chi$ of
integral one in those coordinates and put
\[
 \chi_R(y)=R^{-Q}\chi(\sigma_R(y)).
\]
Let $C$ be a bounded measurable fundamental domain for the
left action of $\Gamma$ on $G$, with translates disjoint
up to null sets. Define
\[
 p_\gamma(x)=\int_{\gamma C}\chi_R(x^{-1}y)\,dy.
\]
These smooth functions form a locally finite,
$\Gamma$-equivariant partition of unity. The values and
derivatives that can contribute at $x$ have their labels
$\gamma$ at distance $O(R)$ from $x$.

Let $X_i$ be the left-invariant field whose value at the
identity is $e_i$. Differentiation in $x$ acts on the kernel
by minus the corresponding right-invariant field in the
variable $x^{-1}y$. Under the transported group law this
right-invariant field is $R^{-w_i}$ times the right-invariant
field generated by $e_i$. The latter fields are uniformly
bounded on the support of $\chi$, since the transported laws
converge smoothly there. Integration over the disjoint
fundamental cells, using the Jacobian $R^Q$, now gives
\begin{equation}\label{eq:deriv}
 \sum_{\gamma\in\Gamma}|X_i p_\gamma(x)|\le C R^{-w_i}.
\end{equation}

The barycentric map $x\mapsto(p_\gamma(x))_\gamma$ takes
$G$ equivariantly into the simplex on the vertex set
$\Gamma$. This simplex is a contractible free simplicial
model for the group, with the left action on vertices. The
pullback of the Whitney form of $\zeta$ \cite{Whitney} is
\begin{equation}\label{eq:whitney}
 \omega_R=\sum_{\gamma_0,\ldots,\gamma_8}
 \zeta(\gamma_0,\ldots,\gamma_8)\,
 p_{\gamma_0}\,dp_{\gamma_1}\wedge\cdots\wedge dp_{\gamma_8}.
\end{equation}
For clarity, the barycentric formula commutes with
differentiation by $\sum_\gamma p_\gamma=1$ and
$\sum_\gamma dp_\gamma=0$. Its integral on an oriented
simplex in the universal model is the value of the
alternating cochain on that simplex. The equivariant
barycentric map therefore pulls this class back to the
class of $\zeta$ on $\Gamma\backslash G$. One can also
compare two partition choices directly by interpolating
them on $G\times[0,1]$.

Let $M_\zeta(s)$ be the supremum on tuples of diameter
at most $s$. Evaluating \eqref{eq:whitney} on
$X_0,\ldots,X_7$, applying \eqref{eq:deriv} in each term
of the wedge expansion, and summing first in the
differentiated labels gives
\[
 \sup_x|\omega_R(X_0,\ldots,X_7)(x)|
 \le C M_\zeta(CR)\prod_{i=0}^7 R^{-w_i}=o(1).
\]
The form descends to the compact nilmanifold. Its integral
over the fundamental class is independent of $R$ and tends
to zero, so its top-degree cohomology class vanishes.
\end{proof}

\begin{corollary}\label{cor:growth}
A class with a cocycle representative of growth $O(R^p)$
lies in $W_p$. A class in $W_p$ with a representative of
growth $o(R^p)$ lies in $W_{p-1}$.
\end{corollary}

\begin{proof}
If $c$ has leading weight $p'>p$, choose by
Proposition~\ref{prop:alg} a class $e\in W_{Q-p'}$ with
$c\smile e\ne0$. Lemma~\ref{lem:integration} represents $e$
by a cocycle of growth $O(R^{Q-p'})$. The homogeneous cup
product with the given representative of $c$ has growth
$O(R^{Q-p'+p})=o(R^Q)$, contradicting Lemma~\ref{lem:top}.
If $c\in W_p$ and its representative has growth $o(R^p)$,
the same argument shows $c\smile W_{Q-p}=0$. Perfectness
on complementary weights gives $c\in W_{p-1}$.
\end{proof}

\section{The associated graded of the induced map}\label{sec:induction}

Suppose $\Psi:\Gamma\to\Gamma$ is a bi-Lipschitz bijection
onto a subgroup $\Gamma'$ of index $k$, regarded as a
quasi-isometry of $\Gamma$ to itself. Write
\[
 \Psi^{c}(x)=\Psi(c)^{-1}\Psi(cx)\qquad(c,x\in\Gamma),
\]
and denote the pointwise closure of these normalized translates
by $\mathscr X_\Psi$. Its elements have common quasi-isometry
constants and fix the identity. Every translate is a bijection
onto $\Gamma'$. This property persists in the closure, since
the common lower quasi-isometry bound puts the preimages of
each fixed point in a finite ball.

\subsection{The dilation parameter}

The rescaled counting measures on $\Gamma$ and $\Gamma'$
converge to Haar measure on $G_\infty$ divided by the
respective covolumes. For example, the rescaled translates of
a bounded fundamental cell have volume $R^{-Q}$ times its
original volume. Their diameters on compact coordinate sets
tend to zero, so the convergence follows by Riemann sums.
The covolume of $\Gamma'$ is $k$ times that of $\Gamma$.

A pointed blow-down $F$ of $\Psi$ is obtained from normalized
translates by rescaling the domain and target, allowing the
centers to vary. These blow-downs are bi-Lipschitz homeomorphisms
of $G_\infty$ by Pansu's results \cite{Pan83,Pan89} and the common quasi-isometry
bounds. For a compactly supported continuous test function
$\varphi$, bijectivity of each normalized translate gives
\[
 R^{-Q}\sum_{x\in\Gamma}\varphi\big(\sigma_R(\Psi^c(x))\big)
 =R^{-Q}\sum_{y\in\Gamma'}\varphi(\sigma_R(y)).
\]
Along a blow-down sequence the left side converges to the
integral of $\varphi\circ F$ divided by the covolume of
$\Gamma$. The common lower bounds localize the summands to
a compact set in rescaled coordinates. The right side gives
the integral of $\varphi$ divided by the covolume of
$\Gamma'$. Thus, writing $m$ for Haar measure on $G_\infty$,
we have $F_*m=k^{-1}m$. At every Pansu differentiability point,
the same measure identity passes to the differential by
translating and rescaling. In particular
\[
 |\det D_x F|=k.
\]
By Proposition~\ref{prop:facts}, the differential has the form
$\delta_t g$ with $t=k^{1/Q}>0$ and
$g\in\{1,\delta_{-1},\tau,\delta_{-1}\tau\}$.

\subsection{An extreme first-layer map}

For a quasi-isometry, Mukherjee, Samanta and Thandar
\cite[Theorem~5.6]{MST} identify the set $\cA(\Psi)$ of
invariant-mean first-layer displacement maps with the convex
hull of the first layers of Pansu differentials of pointed
blow-downs. Each extreme point is such a first layer.
Here the possible differentials belong to the set of four
maps just described. Choose an extreme point $L$ and write
$\widehat L=\delta_t g$ for its graded extension. This
extension is unique because the graded algebra is generated
by its first layer.

By \cite[Theorem~6.3]{MST}, there is an ergodic invariant
probability measure $\mu$ on $\mathscr X_\Psi$ with average
first-layer map $L$, for the action
\[
 (f\cdot x)(y)=f(x)^{-1}f(xy).
\]
The cohomology-ring isomorphism $\mathcal I_L$ from that
theorem is induced by the cochain map
\begin{equation}\label{eq:induction}
 \kappa(\xi)(x_0,\ldots,x_j)=
 \int_{\mathscr X_\Psi}\xi(f(x_0),\ldots,f(x_j))\,d\mu(f).
\end{equation}
Indeed, the coupling of \cite[Prop.~6.2]{MST} has
$\omega_{\Gamma}(x^{-1},f)^{-1}=f(x)$. Substituting this into
Gotfredsen--Kyed's formula \cite[Section~4]{GK} gives
\eqref{eq:induction} in every degree.

We claim that for $\mu$-almost every $f$, the rescalings
of $f$ at the identity converge locally uniformly to
$\widehat L$. The pointwise ergodic theorem \cite{Lin},
applied to the finitely many displacement functions
defining the first-layer map, gives convergence of their
averages over word balls to the values of $L$.
Word balls in the nilpotent group are a tempered
F\o lner sequence by polynomial volume asymptotics.
Take any subsequential blow-down $F$ of such an $f$.
For any fixed bounded rescaled ball, a sufficiently close
pointwise approximation of $f$ by a normalized translate
of $\Psi$ agrees on the whole corresponding finite word
ball. Taking balls of rescaled radius tending to infinity
shows that $F$ is a pointed blow-down of $\Psi$ as well.
The first-layer differential of $F$ lies in $\cA(\Psi)$
almost everywhere. The blow-down averaging identity
\cite[Lemma~5.4]{MST}, applied to $f$ at radii $aR$ for
each positive integer $a$, says that the average of this
differential over the limiting ball of radius $a$ equals
$L$. Since $L$ is extreme, the differential equals $L$
almost everywhere on each such ball. A map with this
constant first-layer differential is the graded
extension $\widehat L$ by \cite[Lemma~4.7]{MST}.
Thus all subsequential blow-downs agree, proving the claim.

\begin{proposition}\label{prop:induction}
The map $\mathcal I_L$ preserves $W$ and satisfies
$\Gr^W\mathcal I_L=\widehat L^*$ under the identification of
Proposition~\ref{prop:alg}.
\end{proposition}

\begin{proof}
For $\omega\in F_p$ closed, Lemma~\ref{lem:integration}
and the common quasi-isometry bounds show that
$\kappa(\xi_\omega)$ has growth $O(R^p)$. Corollary~\ref{cor:growth}
gives preservation of $W_p$.

By Proposition~\ref{prop:alg}, $\widehat L^*\omega^{(p)}$
is the leading part of a closed $\omega'\in F_p$.
To compare cochains, translate the first vertex to the
identity. For any sequence of tuples of diameter at most
$R\to\infty$, pass to a subsequence on which the rescaled
vertices converge to $z_i$. The almost-sure blow-down claim
and Lemma~\ref{lem:integration} show that
\[
 R^{-p}\kappa(\xi_\omega)(x_0,\ldots,x_j)
 \longrightarrow
 \int_{\sigma_\infty(\widehat Lz_0,\ldots,\widehat Lz_j)}
       \omega^{(p)}.
\]
Here dominated convergence applies since the rescaled
integrands are uniformly bounded by the common upper
quasi-isometry bound and the growth of $\xi_\omega$.
The expression on the right equals the corresponding
limit for $\xi_{\omega'}$ because graded automorphisms
commute with exponential coning. The subsequence argument
proves that the difference of these cochains has growth
$o(R^p)$. Its class is in $W_{p-1}$ by
Corollary~\ref{cor:growth}, proving the associated graded
identity.
\end{proof}

\begin{proof}[Proof of Theorem~\ref{thm:main}]
The proposed bijection induces a filtered automorphism of the
cohomology ring. Its associated graded is $(\delta_t g)^*$ and $t^{23}=k$. By
Corollary~\ref{cor:rigid}, $t=1$. Thus $k=1$, contrary to
the subgroup being proper.
\end{proof}

\subsection*{Acknowledgment}
AI assistance was used in preparing the draft and in the
symbolic and independent exact arithmetic checks.

\end{document}